\documentclass{article}
\usepackage[a4paper, margin=2.5cm]{geometry}

\usepackage[utf8]{inputenc}
\newcommand{\dint}{\displaystyle\int}
\newcommand{\dsum}{\displaystyle\sum}
\usepackage{amsmath,amsfonts,amssymb}
\usepackage{stmaryrd}
\usepackage{natbib}
\usepackage{float}
\usepackage{graphicx}
\usepackage{url}
\usepackage[dvipsnames]{xcolor}

\title{Uncovering Data Across Continua: An Introduction to Functional Data Analysis.}

\date{}

\begin{document}

\maketitle

\begin{center}
\textbf{Sophie Dabo-Niang$^{1}$, Camille Frévent$^{2}$}

\vspace{0.1cm}
$ ^{1}$\begin{small}{Univ. Lille, CNRS, UMR 8524-Laboratoire Paul Painlevé, Inria-MODAL, Lille, F-59000, France}\end{small}\\
$^{2}$\begin{small}{Univ. Lille, CHU Lille, ULR 2694 - METRICS: Evaluation des technologies de sant\'e et des pratiques m\'edicales, F-59000 Lille, France}\end{small} \\
\end{center}

\section*{Abstract}
In a world increasingly awash with data, the need to extract meaningful insights from data has never been more crucial. Functional Data Analysis (FDA) goes beyond traditional data points, treating data as dynamic, continuous functions, capturing ever-changing phenomena nuances. This article introduces FDA, merging statistics with real-world complexity, ideal for those with mathematical skills but no FDA background.

\section{Introduction}
Nowadays, advancements in data collection technologies like sensors, computer vision, medical imaging, IoT and wearables have generated vast volumes of high-frequency data across various fields.
These data are not just a collection of numbers and tables but a rich, dynamic tapestry of information that captures the essence of change over a continuum. Functional Data Analysis (FDA)
\citep{bosq2000linear,ramsay2005,ferraty2006nonparametric,kokoszka2017introduction} efficiently handles large-scale, high-dimensional datasets, extracting valuable insights from data containing structured
information.

Unlike traditional statistics dealing with discrete data points, FDA focuses on entire functions, curves or shapes, providing insights into continuous changes.
Whether analyzing time series, spatial data, growth curves, or any structured dataset, FDA excels at capturing ongoing change.
FDA's applications span various fields like medicine, biology, chemistry, economics, and environmental science, offering insights beyond isolated measurements. It aids in patient health tracking, economic trend analysis, chemical or environmental management by modeling and understanding complex systems. In manufacturing, FDA can be applied to monitor continuous processes, such as chemical reactions, quality control measurements, and equipment performance. It helps detect deviations from the desired process behavior \citep{fabioetal20}.
In computational biology, FDA involves studying complex biomolecular structures and understanding the relationship between organism shapes and functionality.
FDA techniques are also applied to analyze longitudinal patient data \citep{yao2005functional}, which are common in clinical trials. This enables the study of disease progression, treatment effectiveness, and personalized medicine.
Furthermore, in biology, in the omics data context, gene expression data comprise measurements of gene expression levels across thousands of genes at multiple time points. 
By considering these data as a function of time, FDA \citep{leng2006classification, cremona2019functional} can help researchers to better understand the general features and dynamics of gene expression, to identify key genes associated with specific diseases or biological processes, and to identify differences or similarities between genes. \\
In economics, FDA is employed to analyze longitudinal data, such as stock prices, gross domestic product trends, and inflation rates. It helps identify long-term patterns, cyclic behavior, and structural changes \citep{horvath2012inference}.
In environmental science, FDA is used to analyze temporal or space-time environmental data, such as temperature records, precipitation patterns, and ocean currents. It aids in understanding long-term climate trends and variability. FDA can be applied to study spatial data, helping to identify pollution hotspots \citep{frevent2023investigating}, or study vegetation growth, and monitor land use changes over time.\\

In essence, FDA transcends traditional data analysis limitations by leveraging data with functionality, providing valuable statistical tools for researchers and professionals seeking deeper insights and solutions to complex problems. We explore FDA's significance, mathematical foundations, practical applications, and future prospects to unveil its transformative potential.

\section{The Significance of Functional Data Analysis}
In various fields, FDA provides a powerful set of methods to model, analyze, and interpret data that exhibit continuous variation, allowing researchers and professionals to gain deeper insights, make more accurate predictions, and informed decisions based on the inherent functional nature of the data. This versatility makes FDA a valuable approach in a wide range of scientific and practical applications \citep{silverman2002applied}.\\ 
Employing mathematical domains like linear algebra, functional analysis, probability and statistics, FDA manipulates and analyzes functions by representing data as observations of random variables in a function space. This allows operations like differentiation, integration, and smoothing, facilitating exploration of data structure and variations.
By treating data as functions, FDA helps uncover hidden patterns, relationships, and trends that would be challenging to discern using traditional statistical methods, leading to more informed decision-making and a deeper understanding of complex phenomena.

\subsection*{FDA versus Multivariate statistics}

\begin{figure}[h!]
\begin{minipage}[t]{0.25\linewidth}
(a) \\
\includegraphics[height=4.5cm]{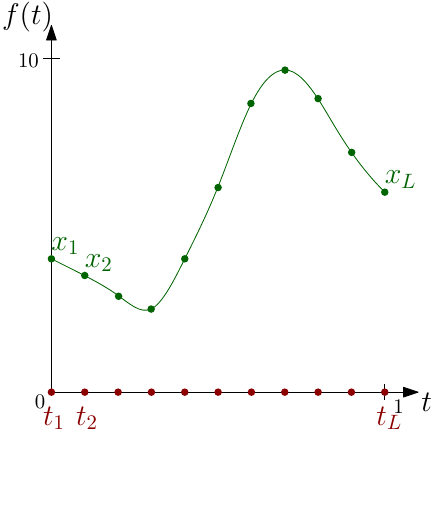}
\end{minipage} \hfill
\begin{minipage}[t]{0.75\linewidth}
(b) \\
\includegraphics[height=4.5cm]{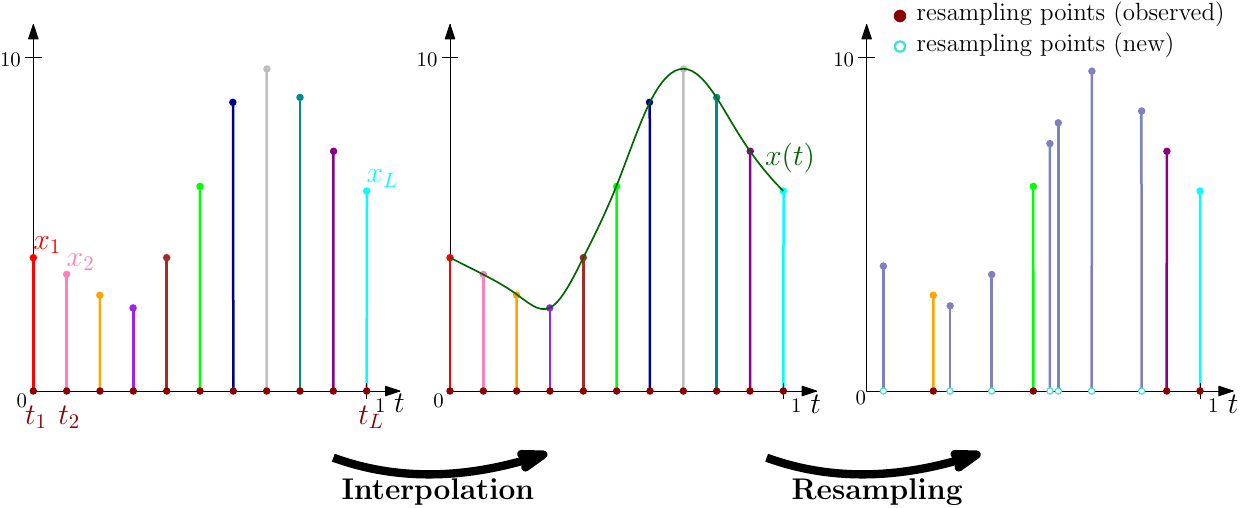}
\end{minipage}
\caption{Example of discrete data from a function $f$}
\label{fig:fda2resamp}
\end{figure}

Is it worthwhile to employ continuous representations, or are we unnecessarily adding complexity to our tasks? Given that discrete data are often needed for computational purposes, what are the benefits of utilizing continuous representations in our analyses?
While discrete data may offer computational convenience, the advantages of working with continuous representations are numerous. By viewing objects as functions, curves, or surfaces, scientists can unlock more powerful analysis techniques, yielding better practical results and more natural solutions. Grenander's principle of discretizing as late as possible underscores the importance of retaining continuous representations for as long as feasible, highlighting their inherent value in data analysis workflows. With this in mind, let us consider how continuous representations enhance our understanding and analysis of data. \\
If data are sampled from an underlying function (e.g., Figure \ref{fig:fda2resamp} (a)), and time points are synchronized across observations, focusing solely on heights, then analysis can be conducted using the vector 
 $\mathbf{x}=(x_1, x_2, \dots, x_L)^\top $. 
If the time points hold significance as well, then it is necessary to retain them alongside the height data: 
$ ((t_1,x_1), (t_2,x_2), \dots, (t_L,x_L))^\top $.
With continuous functions one can interpolate and resample at arbitrary points (e.g., Figure \ref{fig:fda2resamp} (b)) and compare easily observations with different time points, as elements of a function space.

\begin{figure}[H]
\begin{minipage}[t]{0.24\linewidth}
(a) \\
\includegraphics[width=\textwidth]{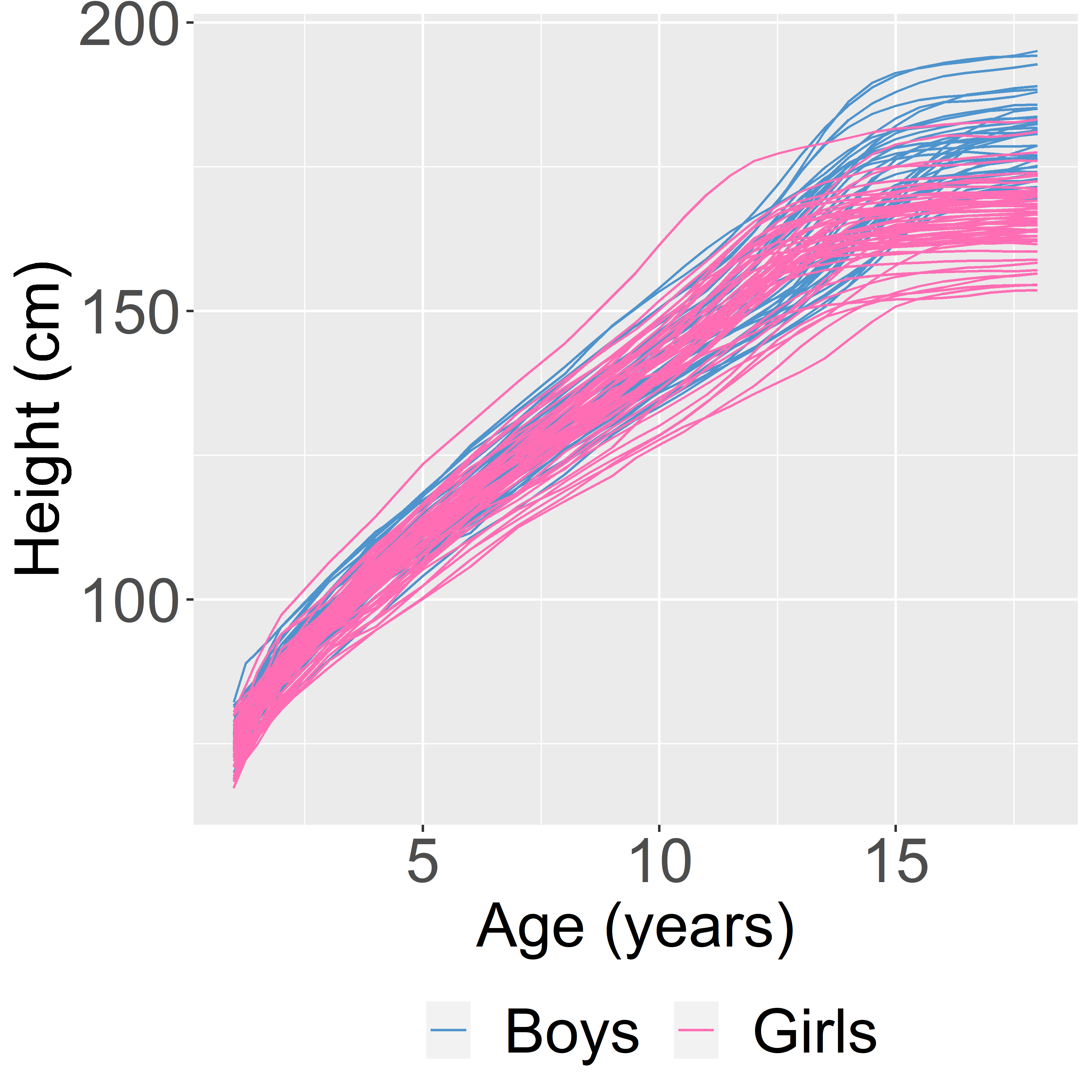}
\end{minipage} \hfill
\begin{minipage}[t]{0.24\linewidth}
(b) \\
\includegraphics[width=\textwidth]{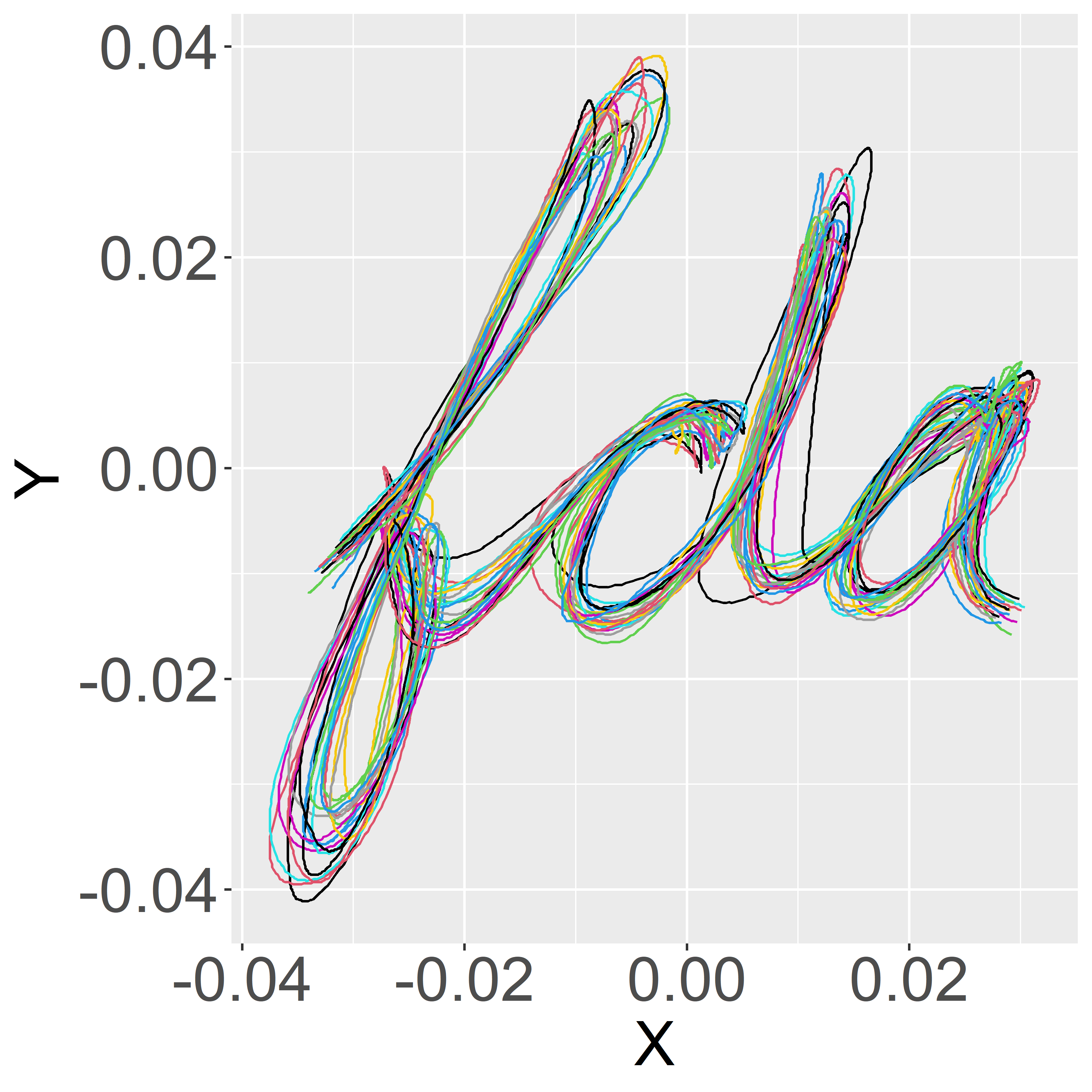}
\end{minipage}
\begin{minipage}[t]{0.24\linewidth}
(c) \\
\includegraphics[width=\textwidth]{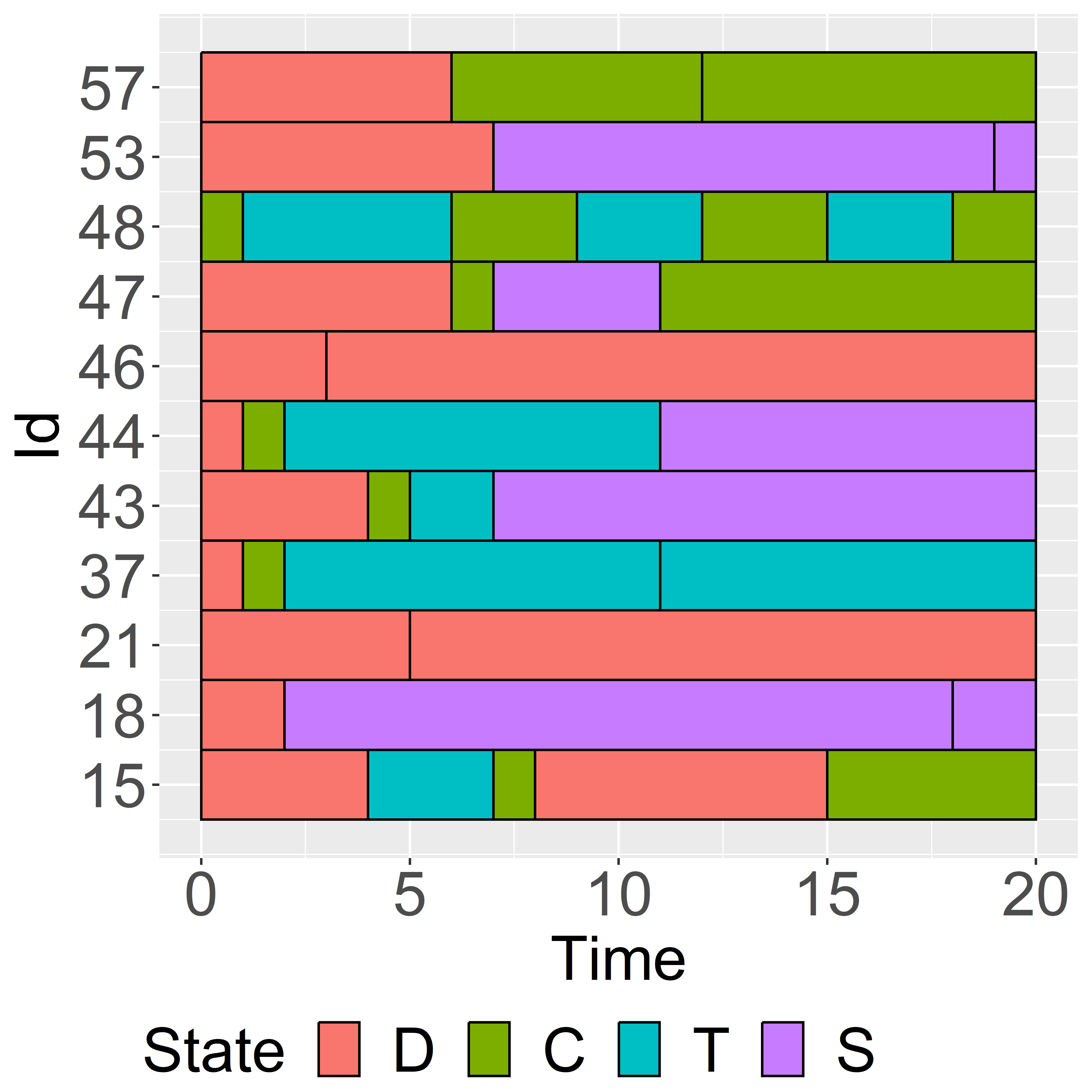}
\end{minipage} \hfill
\begin{minipage}[t]{0.24\linewidth}
(d) \\
\includegraphics[width=\textwidth]{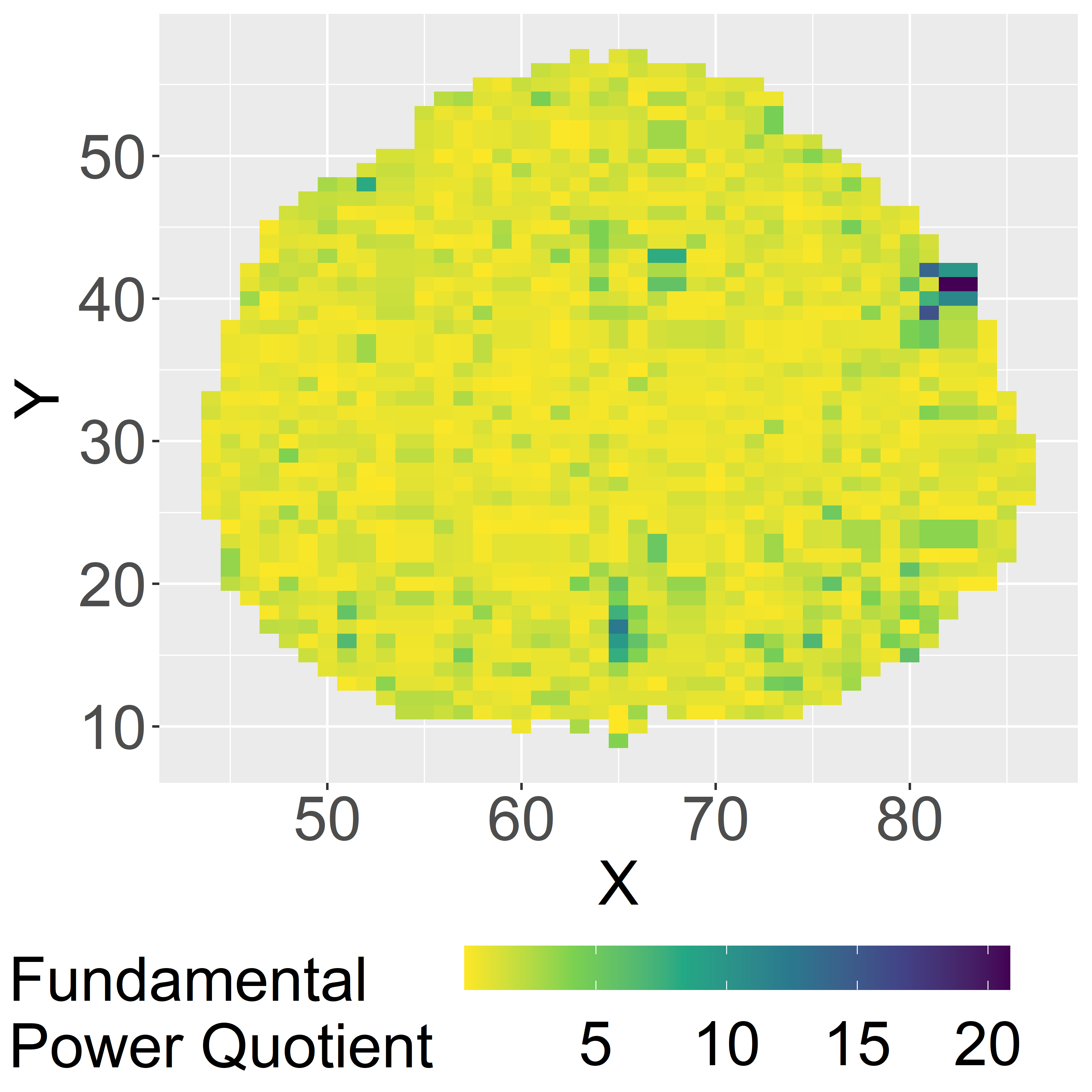}
\end{minipage}
\caption{Examples of functional data. Panel (a) represents the most common case of functional data over time with the growth dataset from the \textsf{R} package \textit{fda}, panel (b) represents the horizontal ($X$) and vertical ($Y$) positions of a pen while writing the word ``fda'' with the handwrit dataset from the \textsf{R} package \textit{fda}, panel (c) represents categorical functional data with the care trajectories of patients over time (care dataset from the \textsf{R} package \textit{fda}) and panel (d) represents functional data in the case of images, here functional magnetic resonance imaging measurements for a human brain with the brain dataset from the \textsf{R} package \textit{gamair}.}
\label{fig:examples}
\end{figure}

In traditional data analysis, one might work with data points in a table where each row represents an observation (e.g., $\mathbf{x}$) and each column represents a variable. In FDA, the data are treated as functions, where each observation is considered as a function (e.g., in Figure \ref{fig:examples} (a)) that maps a continuous variable (often time, frequency, wavelength or a spatial dimension) to a measured value. These functions represent how the data change over the continuum.\\ Understanding functions necessitates a profound grasp of the structures lying beneath them. Analyzing these structures requires a solid foundation of mathematical representations. The FDA approach empowers researchers to investigate various models extensively, thus expanding the comprehension of data characterized by functional structures across diverse fields of science and engineering. Examples of functional data are illustrated in Figure \ref{fig:examples}.

\section{Mathematical Foundations of FDA}
FDA involves a variety of specialized statistical techniques for handling functional data \citep{ramsay2005}, including methods for function smoothing, visualization (e.g., plotting entire functions as curves or surfaces), dimension reduction (e.g., functional principal component analysis), functional regression, and clustering. These techniques account for the continuous nature of the data and are designed to capture underlying patterns and structures in functions. As said before, analyzing these functions involves mathematical representations.

Let $(\Omega, \mathcal{A}, \mathbb{P})$ be a probability space, $\mathcal{F}$ a function space (e.g., a separable Banach space or a Hilbert space). A functional random variable is a variable 
$$
X=\{X(t), t\in\mathcal{T}\} : \Omega \to \mathcal{F},
$$
taking values in $\mathcal{F}$ (of eventually infinite dimension). A functional data is then an observation of the functional random variable $X$.\\

If $\mathcal{T} \subseteq \mathbb{R}$ then $X$ is a curve while an image may be considered as a functional data in the case where $\mathcal{T} \subseteq \mathbb{R}^2$. If $\mathcal{T} \subseteq \mathbb{R}^d$ ($d>2$), $X$ has a more complex structure.

Let us consider in the following the commonly used functional Hilbert space $\mathcal{L}^2(\mathcal{T},\mathbb{R})$, the space of $p$-dimensional vector-valued square-integrable functions on $\mathcal{T}$, and give the main background to analyse functional data. First, consider the inner product on that Hilbert space: for $f,g \in \mathcal{L}^2(\mathcal{T},\mathbb{R})$, 
$\langle f,g \rangle = \int_{\mathcal{T}} f(t) \ g(t) \ \text{d}t$. The mean and the covariance functions of the random variable $X \in \mathcal{L}^2(\mathcal{T},\mathbb{R})$, assumed as smooth functions, are respectively 
$$ \mu_X(t) = \mathbb{E}[X(t)] \in \mathbb{R}$$ and 
$$ \text{C}_X(s,t) = \text{Cov}(X_s,X_t) = \mathbb{E}\{[X(t) - \mu_X(t)] [X(s)-\mu_X(s)] \} \in \mathbb{R}.$$
The latter is viewed as the kernel of the linear Hilbert-Schmidt operator $\Gamma_X$ on $\mathcal{F}=\mathcal{L}^2(\mathcal{T},\mathbb{R})$:
$\Gamma_X: \mathcal{F} \rightarrow \mathcal{F},\; 
\Gamma_X f(t)=\int_\mathcal{T} \text{C}_X(t,s) f(s) \ \text{d}s $. Note that $\Gamma_X$ admits the spectral decomposition
$
\Gamma_X=\sum_{j\ge 1} \lambda_{j} f_j\otimes f_j
$
where $(f\otimes g)(x)=\langle f,x \rangle g$, $x\in \mathcal{F}$, $\{f_j\}_j$ is a complete orthonormal system in $\mathcal{L}^2(\mathcal{T},\mathbb{R})$ and 
$\{\lambda_{j}\}_j$ is a decreasing sequence of positive real numbers such that $\sum_{j\ge 1}\lambda_j < \infty$.\\

Let $X_1, \dots, X_n$ be an independent and identically distributed (i.i.d.) sample of $X$.\\
The usual estimator of $\mu_X$ is the method of moments estimator given by
$\widehat{\mu}_X(t) = \frac{1}{n}\sum_{i=1}^n X_i(t).$

In this i.i.d. framework, there are several theoretical guarantees regarding the convergence of $\widehat{\mu}_X$ to $\mu_X$ (such as the law of large numbers, the central limit theorem, and concentration inequalities of the Bernstein type, see Chapter 2 of \cite{bosq2000linear}). For instance, for the Hilbert space $\mathcal F$ equipped with a norm $\|.\|_{\mathcal F}$, if 
$\mathbb E\left[\|X_1\|_{\mathcal F}\right]<\infty$
then $\widehat{\mu}_X \to \mu_X$ almost surely as $n \to \infty$.
If 
$\mathbb E[\|X_1\|^2_{\mathcal F}]<\infty 
$, $\widehat{\mu}_X$ is asymptotically normally distributed: 
 $ \sqrt{n} \left( \widehat{\mu}_X - \mu_X \right)$ converges in distribution to $\mathcal{N}(0,\text{C}_X)$. \\ 
Classic empirical estimators of the covariance operator $\Gamma_X$ and covariance function $C_X$ are
$
\widehat{\Gamma}_X=\frac{1}{n}\sum_{i=1}^{n} (X_i-\widehat{\mu}_X) \otimes (X_i-\widehat{\mu}_X)
$
and
$
\widehat{C}_X(s,t)=\frac{1}{n}\sum_{i=1}^{n} [X_i(t)-\widehat{\mu}_X(t)] [X_i(s)-\widehat{\mu}_X(s)]$.
Several asymptotic results on $\widehat{\Gamma}_X$ are given in Chapter 4 of \cite{bosq2000linear}. More theoretical details are given in this last reference, as well as in \cite{10.1214/009053606000000272} and \cite{10.3150/20-BEJ1209}.

\paragraph{From raw data to functional data:}
Note that, in practice, we observe raw data (e.g., the average daily temperature in spatial locations described by the first panel of Figure \ref{fig:rawliss}: the temperature in each location is measured every day from 1960 to 1994) of the form
$$x_{i,t_{i,l_i}},\quad t_{i,l_i} \in\mathcal{T},\quad i=1,\dots,n \quad l_i=1\dots, L_i.
$$

\begin{figure}[h!]
\centering
\includegraphics[width=0.6\linewidth]{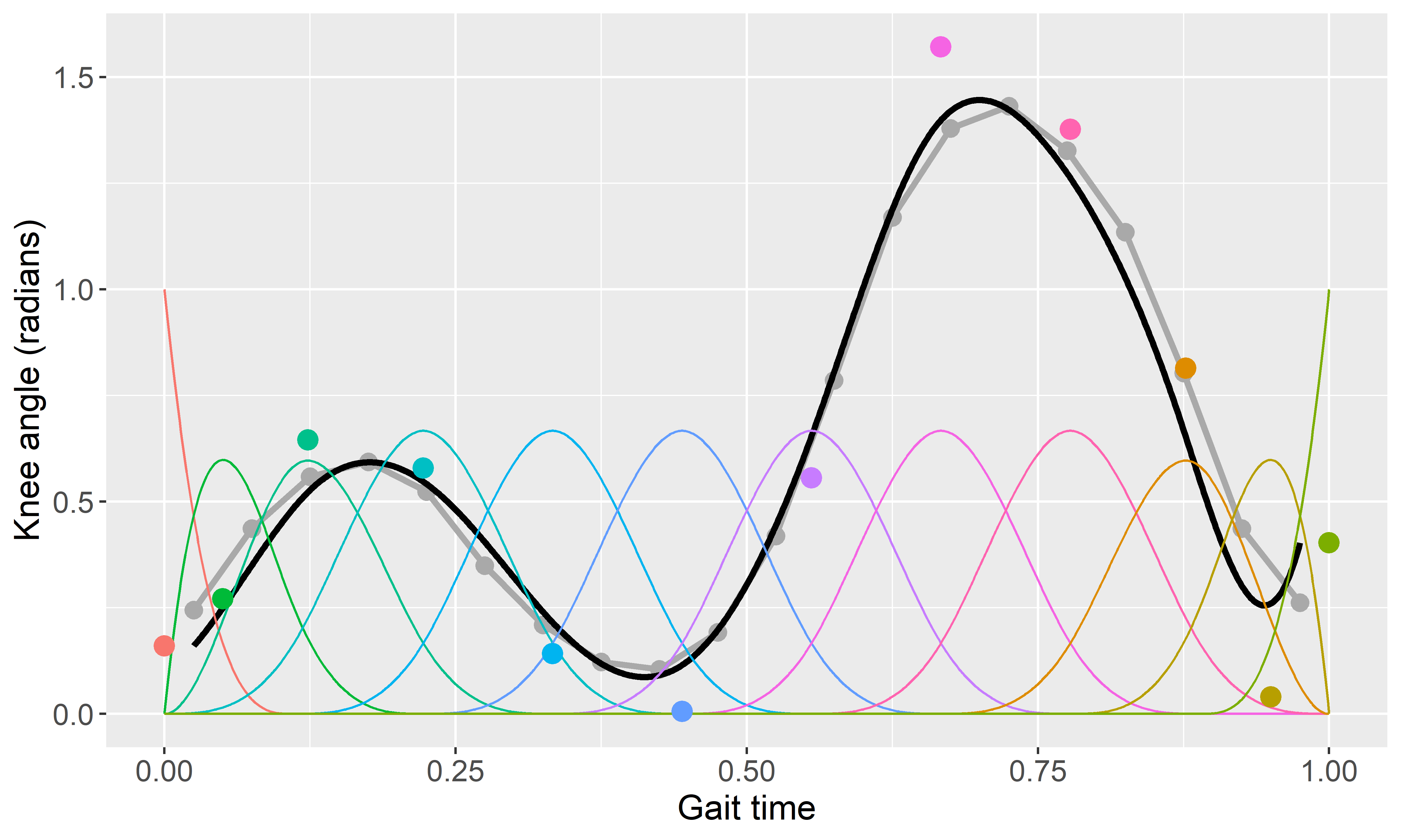}
\caption{Example of transition from raw data (grey connected points) to functional data (black curve) using a cubic $B$-splines basis (colored curves) on the gait dataset from the \textsf{R} package \textit{fda}. The colored points correspond to the $B$-splines coefficients (with the same color as the splines).}
\label{fig:example}
\end{figure}

It should be noted that the observation times $t_{i,l_i}$ can vary in number and value depending on the individual $i$. \\
Following \cite{zhangwang}, we should distinguish the case of sparse and dense functional data. Dense functional data are characterized by the fact that all $L_i$ are larger than some order of $n$. In this case, it is possible to use a smoothing technique on the raw data $x_{i,t_{i,l_i}}$ to recover the original curves. \\
It is common in FDA to assume that the $L_i$ observations $\{x_{i,t_{i,1}}, \dots, x_{i,t_{i,L_i}}\}$ are noisy observations of the smooth latent curve $X_i(.)$. Namely, we have $x_{i,t_{i,l_i}}=X_i(t_{i,l_i})+\varepsilon_{i,t_{i,l_i}}$, where the error terms $\varepsilon_{i,t_{i,l_i}}$ are zero mean and i.i.d. In the early stages of FDA, this smoothing is typically conducted as an initial step by kernel smoothing, local polynomial smoothing, Fourier, spline or penalized spline approaches.\\
The classic smoothing approach is basis expansion by assuming that $X_i, 1 \le i \le n$ can be expressed as a finite combination of the first $K$ functions of a basis functions $\{\phi_{1}, \dots \phi_{K}, \dots\}$ of $\mathcal{L}^2(\mathcal{T},\mathbb{R})$:
$$ X_i(t) = \dsum_{k=1}^{K} a_{i,k} \phi_{k} (t). $$
This is equivalent to write
$ X_i(t) = \phi(t)a_i $
where $\phi(t) = (\phi_{1}(t), \cdots, \phi_{K}(t))$ 
and $a_i = (a_{i,1}, \dots, a_{i,K})^\top$. Then the estimators of the mean and covariance functions of $X$ can be defined respectively on $\mathcal{T}$ by 
$$ \widehat{\mu}_X(t) = \dfrac{1}{n} \sum_{i=1}^n X_i(t)=\dfrac{1}{n} \phi(t) \dsum_{i=1}^n a_i, \text{and } \widehat{\text{C}}_X(s,t) = \dfrac{1}{n} \phi(t) A^\top A \phi(s)^\top
$$
where $A$ is the $n \times K$ matrix whose $i^{\text{th}}$ row $A_i$ is equal to $\left( a_i - \frac{1}{n} \sum_{k=1}^n a_k \right)^\top$. Depending on the nature of the data, various choices for the $\phi_k$ are possible. In the case of periodic data, a Fourier basis is appropriate, whereas for non-periodic data, possible choices are polynomial basis or splines basis. Figure \ref{fig:example} shows a transition from raw data to functional data using a cubic $B$-splines basis. This process is not applicable for sparse functional data due to the very limited quantity of information available for each curve. Sparse data require more sophisticated approaches not discuss here. For more details concerning basis options, please see \cite{ramsay2005} and for other smoothing techniques or more details on the sampling design of functional data, see \cite{zhangwang}. \\

With the growing popularity of functional data analysis, numerous statistical methods have been extended and adapted to this context. In the following discussion, we explore some of the most valuable and widely used ones (principal component analysis, clustering, linear and non-linear regressions). It should be noted that we only consider the case of univariate functional data, i.e., for all $t \in \mathcal{T}, X(t) \in \mathbb{R}$. The structure of multivariate functional data ($X(t) \in \mathbb{R}^p, p \ge 2$) is more complex, please refer to \cite{koner2023second} for more details.

\paragraph{Functional principal component analysis:}

\cite{kleffe1973principal} has expanded upon the conventional and widely employed statistical approach known as Principal Component Analysis (PCA) to accommodate variables with values in a separable Hilbert space. Then PCA has been extended to accommodate specifically univariate functional data (see for example \cite{locantore1999robust} and \cite{ramsay2005}).
This technique is particularly useful for reducing the dimensions of the data, extracting their main features, providing an indication of their complexity and gaining insights into the underlying patterns and structures.
The subsequent paragraph discusses the functional PCA method within the context of univariate functional data, noting that the case of multivariate functional data has also been investigated in the literature \citep{ramsay2005, berrendero2011principal}.

The above spectral decomposition of $\Gamma_X$ is linked to the PCA on the $X_i$. In fact, functional PCA aims to represent the i.i.d. curves $X_i$ using a few ($P$) principal orthogonal eigenfunctions $f_j \in \mathcal{L}^2(\mathcal{T},\mathbb{R})$ so that $$X_i(t) \approx {\mu}_X(t) + \dsum_{j=1}^{P} c_{i,j} f_j(t).$$
This is an approximation of the Kosambi-Karhunen-Loève \citep{kosambi2016decomp,loeve1945functions,karhunen1947under} 
expansion that states that 
$$ X_i(t) = \mu_X(t) + \dsum_{j\ge 1} c_{i,j} f_j(t),$$
where $\{f_j\}_{j\ge 1}$ and $\lambda_{1}\ge \lambda_{2} \ge \dots \ge 0$ are the eigenfunctions and eigenvalues of $\Gamma_X$.
The $c_{i,j}=\langle X_i-{\mu}_X , f_j \rangle$ are called the scores, they extract the main features of $X_i$ and are centered pairwise uncorrelated random variables. \\
Note that the Kosambi-Karhunen-Loève expansion is related to Mercer’s theorem \citep{mercer1909functions}, 
that states
$$
\text{C}_X(s,t)=\dsum_{j\ge 1} \lambda_{j} f_j(t)f_j(s).
$$

Since the mean and covariance functions are unknown, in the early stage of FDA, applying a functional PCA is in practice equivalent to find estimated orthogonal eigenfunctions $\widehat{f}_j$ so that 
$$ \forall t, \dint_{\mathcal{T}} \widehat{\text{C}}_X(s,t) \widehat{f}_j(s) \ \text{d}s = \widehat{\lambda}_j \widehat{f}_j(t). $$
Hence, assuming that $\widehat{f}_j$ can be expressed as
$\widehat{f}_j(t) = \phi(t) b_j$, the task is to find $\widehat{\lambda}_j \in \mathbb{R}$ and $b_j \in \mathbb{R}^K$ so that 
$$ 
\dfrac{1}{n} A^\top A \dint_{\mathcal{T}} \phi(s)^\top \phi(s) \ \text{d}s \ b_j = \widehat{\lambda}_j b_j.
$$
By defining $W = \int_{\mathcal{T}} \phi(s)^\top \phi(s) \ \text{d}s$ and $u_j = b_j^\top W^{1/2}$, and then multiplying the previous equation on the left by $W^{1/2\top}$, we arrive at the classic PCA formula (where $Z = A W^{1/2}$):
$$ 
\dfrac{1}{n} Z^\top Z u_j^\top = \widehat{\lambda}_j u_j^\top.
$$
Subsequently, it becomes straightforward to ascertain the values of $u_j$ and $\widehat{\lambda}_j$. Additionally, the determination of $b_j$ (and consequently $\widehat{f}_j$) can be inferred through the following relationships:
$$ b_j = (u_j W^{-1/2})^\top,\; \widehat{f}_j(t) = \phi(t) b_j.$$
Ultimately, the estimated scores $\widehat{c}_{i,j}$ are given by $\widehat{c}_{i,j} = \langle X_i-\widehat{\mu}_X, \widehat{f}_j \rangle = A_i W b_j.$
It should be noted that, although not discussed here, other approaches for functional PCA have been proposed in the literature. A more complete review can be found in \cite{shang2014survey}.

\paragraph{Functional linear regression:}

Numerous studies have explored regression modeling within the context of functional data \citep{cuevas2002linear,james2002generalized,morris2015functional}.
In the following, we focus on presenting generalized functional linear models, which are designed to model a continuous response variable as a function of functional covariates.\\
In this framework, we consider a real-valued response variable $Y$ and a functional covariate $\{X(t), t \in \mathcal{T} \} \in \mathcal{L}^2(\mathcal{T},\mathbb{R})$.
Throughout the section, we assume that $X$ has been centered and that we have a sample $(Y_i,\{X_i(t), t \in \mathcal{T} \})_{i=1, \dots, n}$ of $n$ i.i.d. replications of $(Y,\{X(t), t \in \mathcal{T} \})$. 

Generalized functional linear regression posits that the relationship between the response variable and the functional covariate is defined as follows:
$$ \mathbb{E}[Y|\{X(t), t \in \mathcal{T} \}] = g^{-1}(\eta),\; \text{Var}[Y|\{X(t), t \in \mathcal{T} \}] = V(g^{-1}(\eta)),$$
where $g$ is a monotonic ``link function'', $V$ is a positive ``variance function'' and $\beta \in \mathcal{L}^2(\mathcal{T},\mathbb{R})$, $\eta$ is a linear predictor defined by 
$ \eta = \alpha + \dint_{\mathcal{T}} X(t) \ \beta(t) \ \text{d}t.$\\
The model finds practical applications in various scenarios, such as establishing associations between the incidence of respiratory diseases (e.g., asthma, lung cancer) and air pollution levels in the months or years leading up to the study. In such cases, a generalized functional linear Poisson regression model is often employed, where the link function $g$ is the logarithm, and $V$ is the identity function.\\
The simplest and most popular model is the so-called functional linear model, where $g$ is the identity function, and $V$ is a constant function:
$$
Y_i= \alpha + \dint_{\mathcal{T}} X_i(t) \ \beta(t) \ \text{d}t +\varepsilon_i.
$$
The random variables $\varepsilon_i$ are assumed i.i.d., scalar variables with a mean of zero and a constant variance. Sometimes, an additional assumption of normality is made.\\
More generally, this linear model may encompass both functional ($X_i$) and non-functional ($Z_i \in \mathbb{R}^d$, $d\ge 1$) covariates, so that:
$$ Y_i = \alpha + Z_i^\top \theta+ \dint_{\mathcal{T}} X_i(t) \ \beta(t) \ \text{d}t + \varepsilon_i,$$
with $\theta\in \mathbb{R}^d.$ \\

The primary challenge in FDA lies in dealing with the infinite dimension of the functional variable. A frequently adopted solution is to approximate $X_{i}$ as a finite combination of orthogonal basis functions (as mentioned earlier), as well as $\beta$.
However, in practice, finding such a basis is not always straightforward. 
An orthonormal basis can be derived through functional PCA:
$$ X_{i}(t) \approx \dsum_{j = 1}^{P} c_{i,j} f_j(t)$$ where the $f_j$ are the orthonormal eigenfunctions. By assuming that $\beta$ can also be written as
$ \beta(t) \approx \sum_{j = 1}^{P} d_{j} f_j(t),$
we then obtain the following truncated linear regression model: $ Y_i = \alpha + Z_i^\top \theta+ \sum_{j = 1}^{P} c_{i,j} \ d_j+ \varepsilon_i.$
Beyond the linear model, this truncation procedure leads to a classic generalized linear model with covariates $c_{i,j}$ and $Z_{i}$.\\

\paragraph{Functional data clustering:}
Clustering is the process of organizing observations into clusters, where observations within each cluster share similar characteristics, while the characteristics of each cluster are distinct from those of others. Clustering methods can be broadly categorized into hierarchical, partitional, and model-based approaches. Researchers have explored adaptations of these categories to the functional data framework. In the case of hierarchical methods, a significant challenge arises in devising an appropriate similarity measure for functional observations.
One approach to addressing this challenge was presented by \cite{hitchcock2007effect}.
Among partitional methods, the most well-known technique is the $K$-means algorithm. It starts by randomly selecting $K$ points as the initial centers of $K$ groups and then assigns each observation to the group with the closest centroid. The centroids of the $K$ groups are then recalculated, and the observations are reassigned to the groups iteratively until convergence.

Then, the $K$-means algorithm relies on measuring the distance between observations, typically using the Euclidean distance for non-functional data. However, when dealing with functional data, this distance metric needs to be adapted. \cite{garcia2015k} conducted a study where they compared various approaches for modifying the $K$-means algorithm in the context of functional data.\\
Model-based clustering based on mixture of distributions have also been proposed. The interested reader may refer to \cite{zhang2023review} for a detailed examination of clustering approaches specific to functional data.\\
For more comprehensive details, methodologies, and applications, please refer to the reviews provided by \cite{horvath2012inference,wang2016functional, koner2023second}.\\

\section{Applications and Future Directions}

A multitude of methods for handling functional data have been introduced, with many others yet to be discovered. In this section, we highlight the potential of functional data through an illustration of clustering using the well-known Canadian Weather dataset from the \textsf{R} package \textit{fda}. We focus on the average daily temperature recorded every day from 1960 to 1994 in 35 spatial locations in Canada.

\begin{figure}[h!]
\begin{minipage}{0.49\linewidth}
\centering
\includegraphics[width=0.8\linewidth]{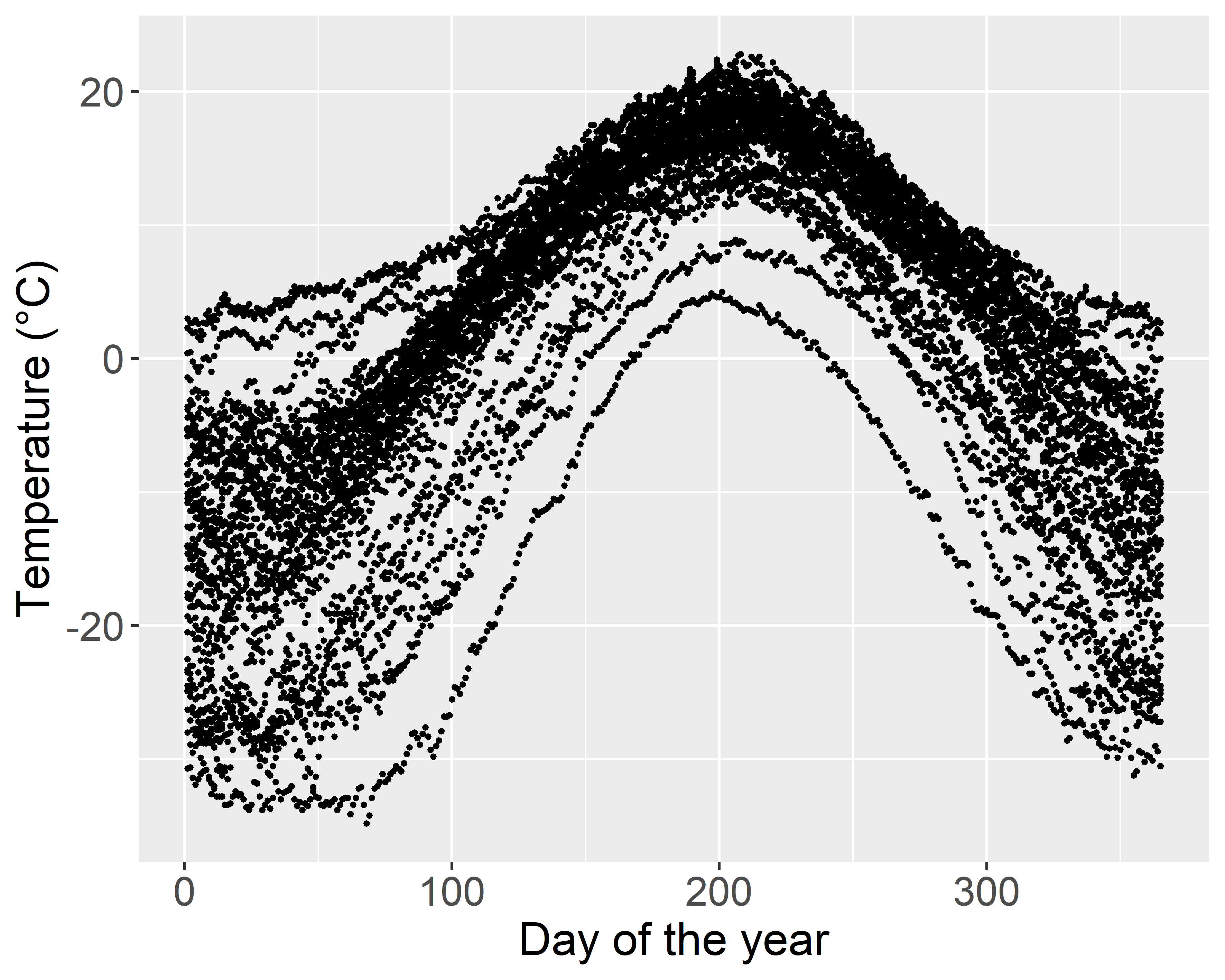}
\end{minipage}
\hfill
\begin{minipage}{0.49\linewidth}
\centering
\includegraphics[width=0.8\linewidth]{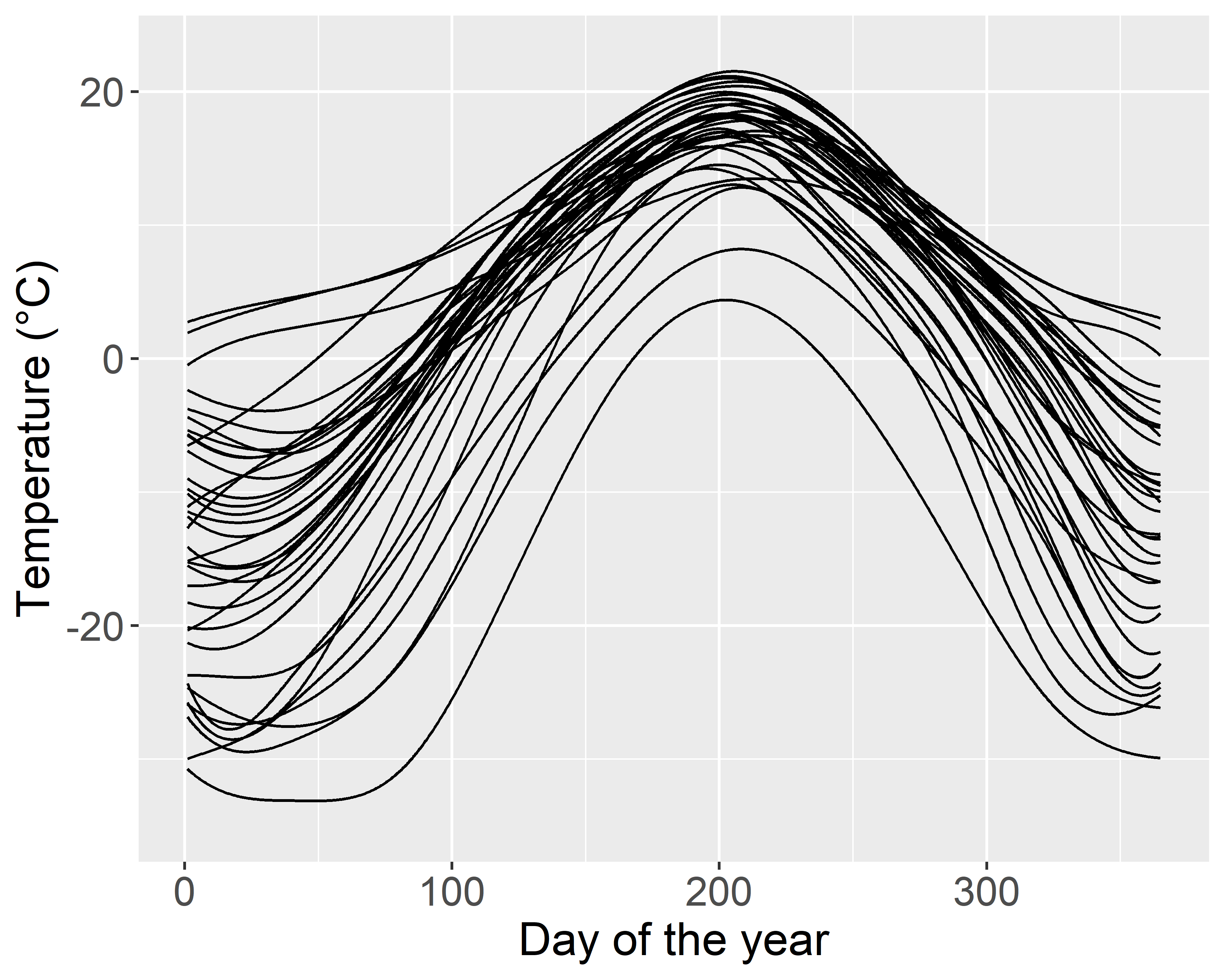}
\end{minipage}
\caption{Observed average daily temperature (left panel) and reconstructed functions (using a $B$-splines basis) for the average daily temperature (right panel) every day from 1960 to 1994 in 35 spatial locations in Canada}
\label{fig:rawliss}
\end{figure}

As said above, in practical applications, functional data are typically observed at discrete points, such as the 365 days of the year. However, it is possible to reconstruct the underlying functions by representing them in a basis of functions. The initial step of this process is depicted in Figure \ref{fig:rawliss}, where a $B$-splines basis has been employed.

\begin{figure}[h!]
\centering
\includegraphics[width=0.4\linewidth]{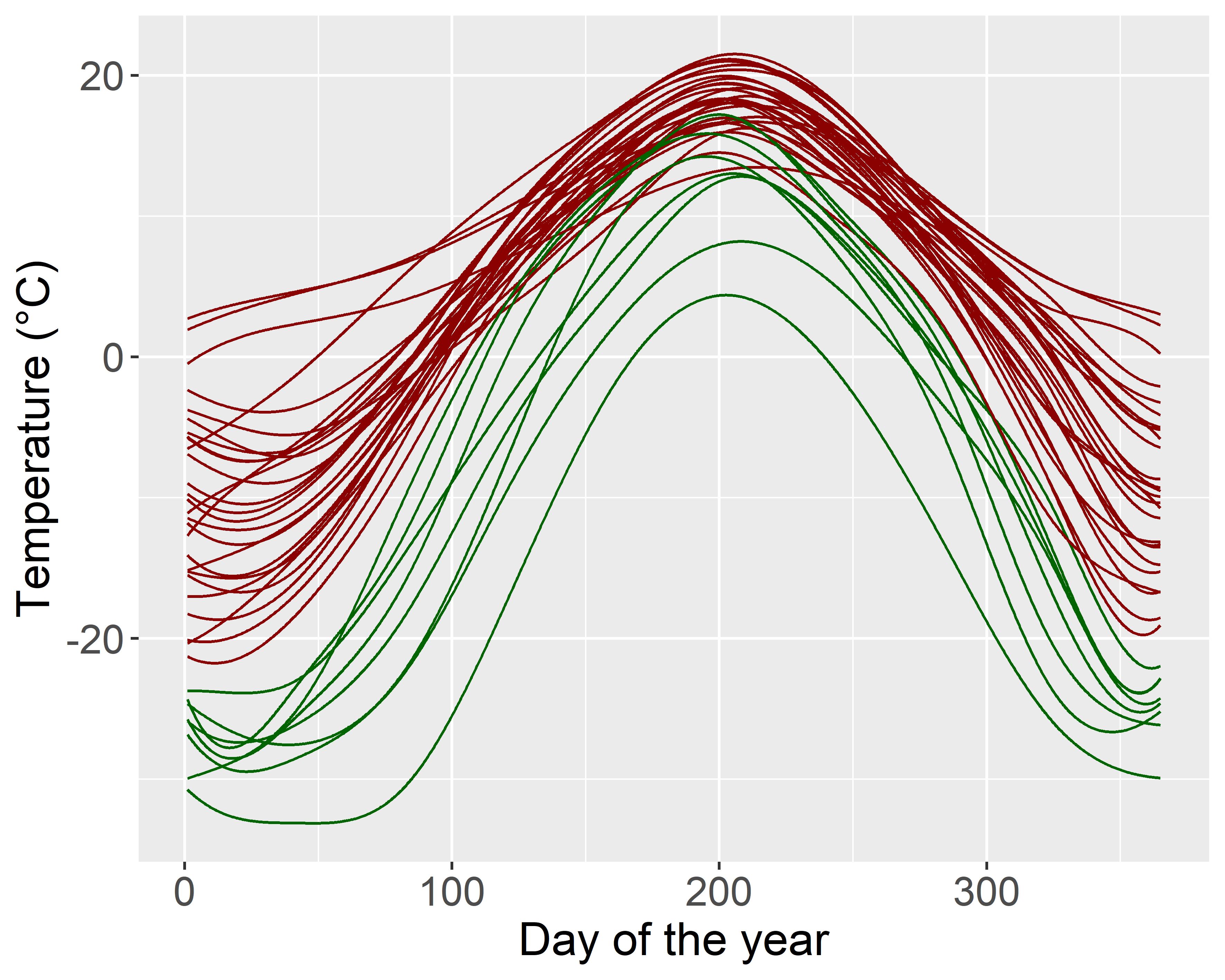}
\caption{Results of the $K$-means algorithm on the Canadian Weather dataset (temperature data) of the \textsf{R} package \textit{fda}}
\label{fig:clustering}
\end{figure}

To distinguish groups of Canadian cities based on their temperature patterns, we applied the $K$-means algorithm proposed by \cite{sangalli2010k} and implemented in the \textsf{R} package \textit{fdacluster}. The optimal number of groups was determined using the optimal average silhouette index and the resulting clustering results are depicted in Figure \ref{fig:clustering}.
We can observe that two distinct groups of Canadian cities emerge from the analysis: the first group (in green) corresponds to cities with consistently lower temperatures throughout the year, while the second group (in red) represents cities with consistently higher temperatures. 

In addition to the functional aspect of the data, the spatial dimension is becoming increasingly relevant, particularly in the context of environmental data. Thus, the literature has seen the emergence of numerous methods specifically tailored to the analysis of spatial functional data. Recently, several spatial cluster detection methods have been introduced in this context. These methods can be used, for instance, to identify environmental hotspots characterized by elevated levels of certain pollutants.

\begin{figure}[h!]
\begin{minipage}[t]{0.64\linewidth}
\begin{minipage}[t]{0.49\linewidth}
\centering
\includegraphics[width=\linewidth]{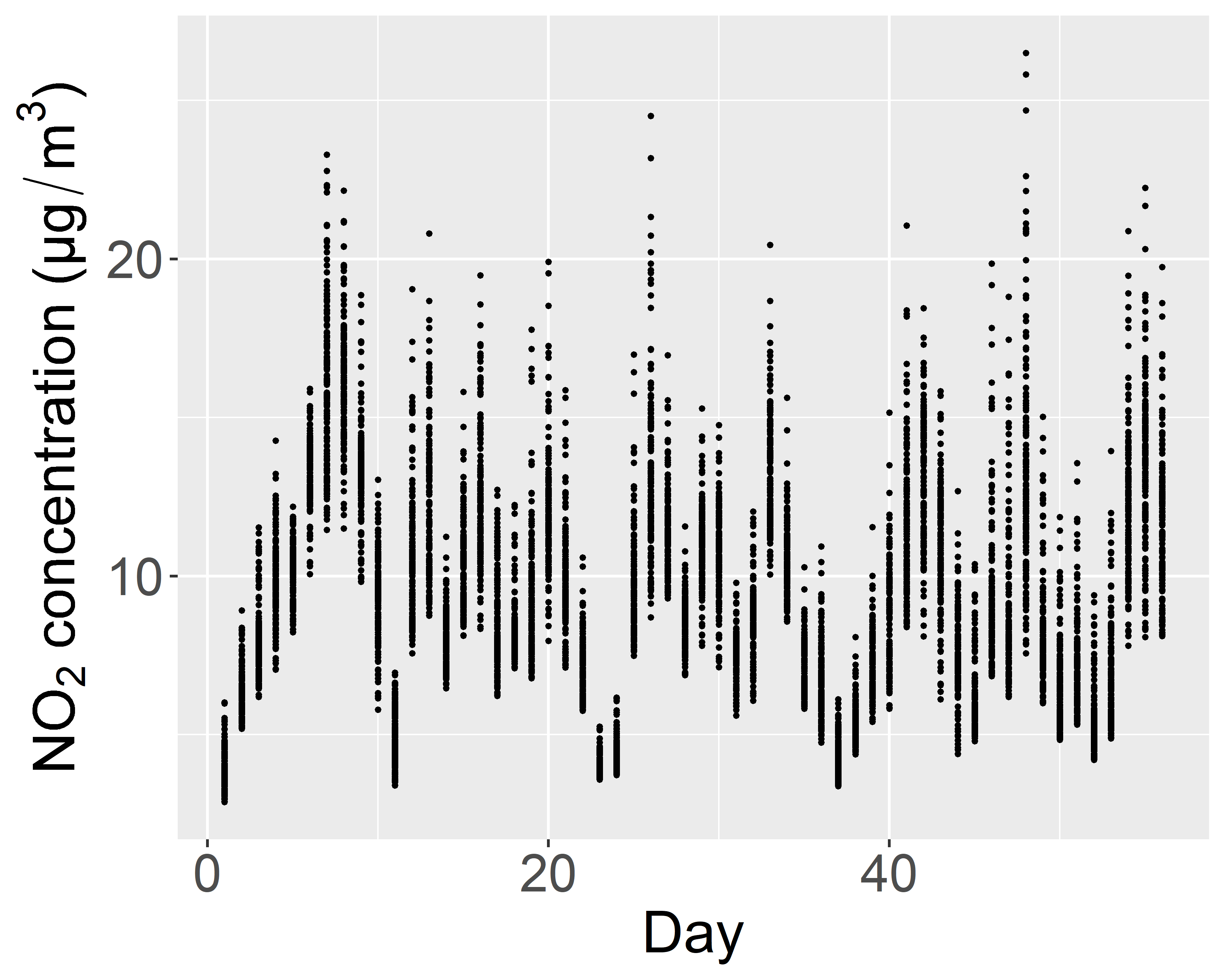}
\end{minipage}
\hfill
\begin{minipage}[t]{0.49\linewidth}
\centering
\includegraphics[width=\linewidth]{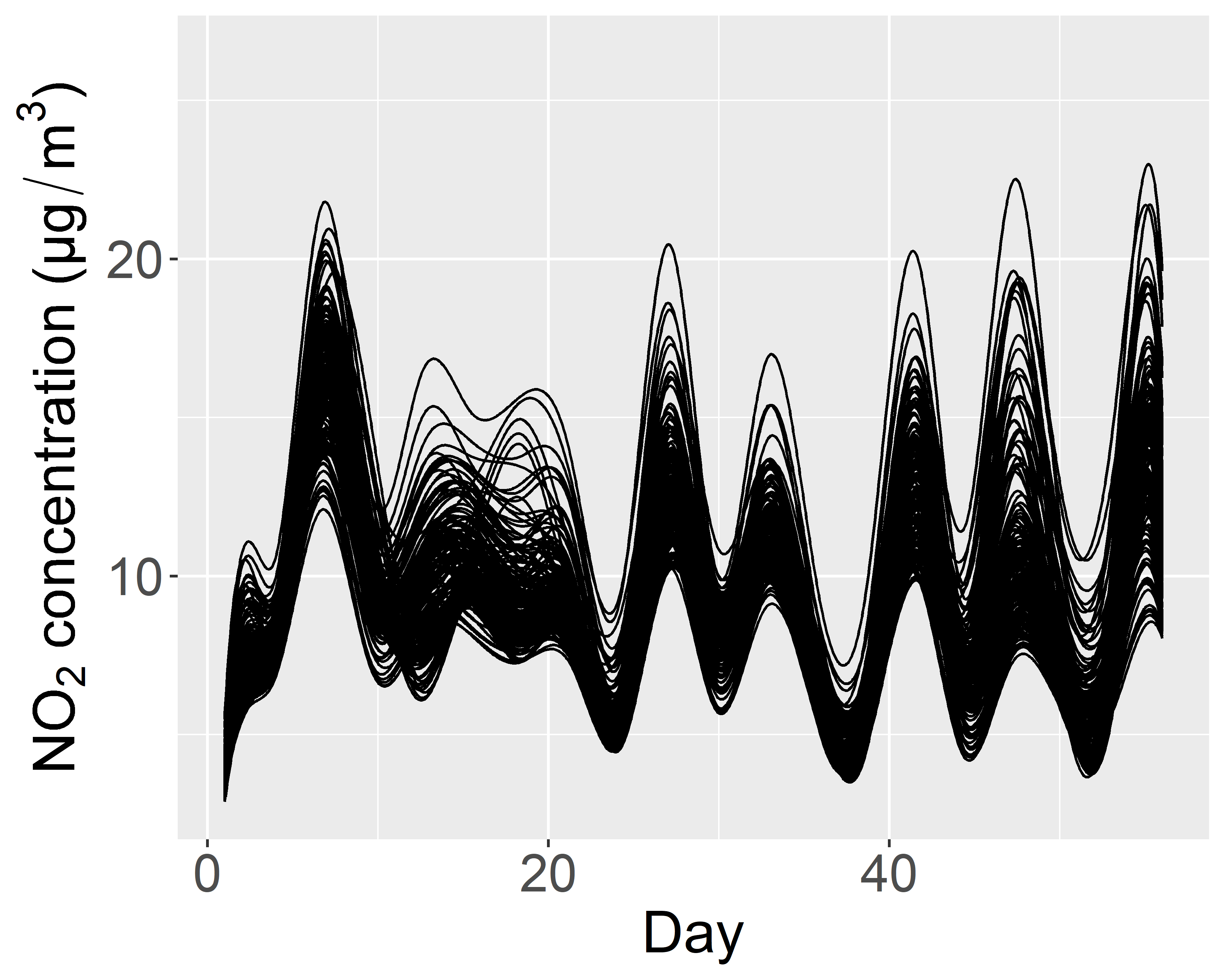}
\end{minipage}
\caption{Observed daily average concentration (left panel) and reconstructed functions (using a $B$-splines basis) for the daily average concentration (right panel) of $\text{NO}_2$ from May 1 to June 25, 2020 in northern France}
\label{fig:rawlisspollution}
\end{minipage}
\hfill
\begin{minipage}[t]{0.32\linewidth}
\centering
\includegraphics[width=\linewidth]{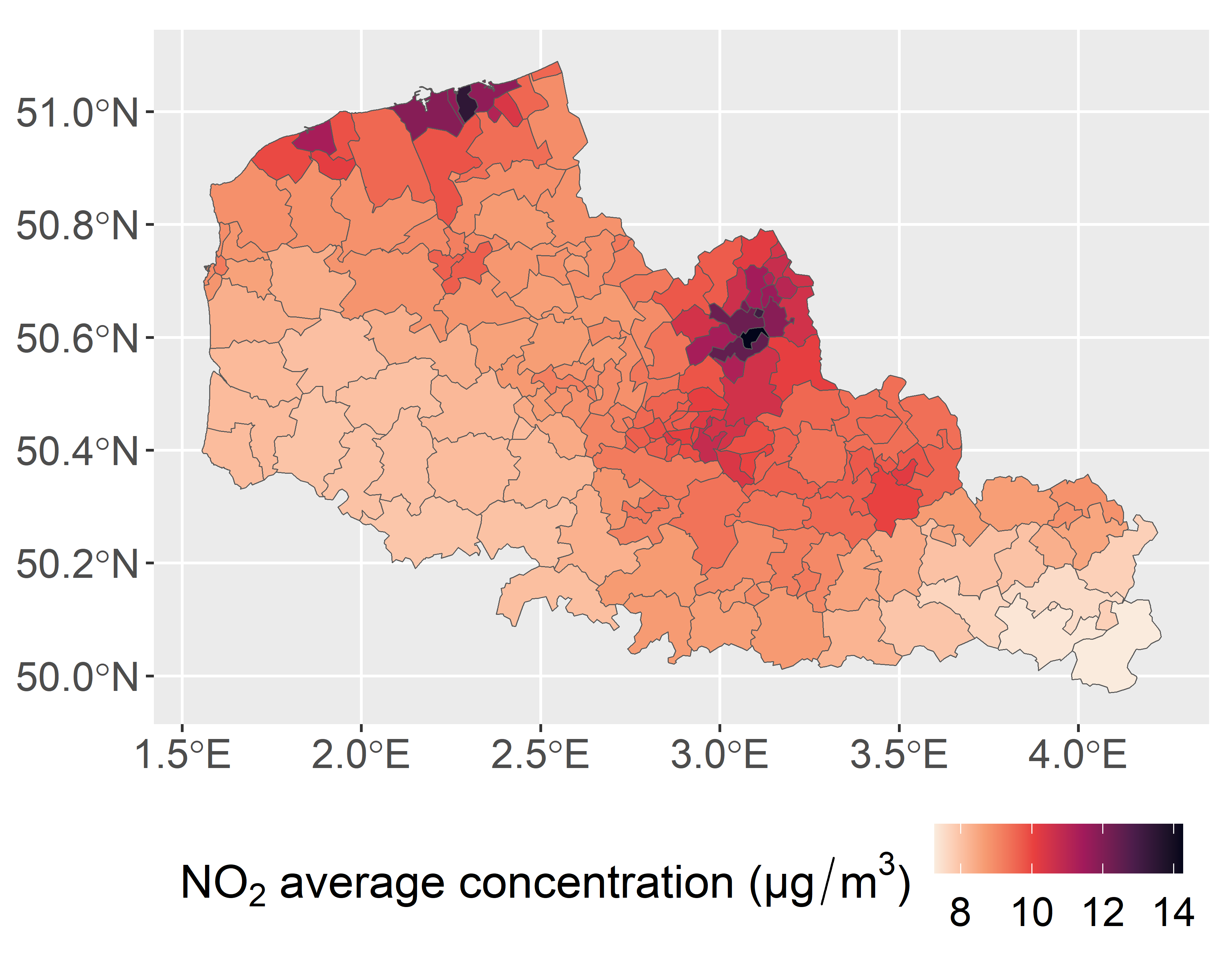}
\caption{Spatial distribution of average $\text{NO}_2$ concentration over the period May 1 to June 25, 2020 in northern France}
\label{fig:mappollution}
\end{minipage}
\end{figure}

In the following example, we will demonstrate a cluster detection approach that incorporates both the spatial and the functional nature of the data. This approach is the distribution-free functional spatial scan statistic (DFFSS) proposed by \cite{frevent2021detecting}, which has been implemented in the \textsf{R} package \textit{HDSpatialScan}. The data used in this example are sourced from the National Air Quality Forecasting Platform (\url{www.prevair.org}) and are available within the package. They comprise the daily average concentration of the pollutant $\text{NO}_2$ recorded from May 1 to June 25, 2020, in northern France. Figures \ref{fig:rawlisspollution} and \ref{fig:mappollution} present the raw data, their functional reconstruction using a $B$-splines basis, and their spatial distribution, respectively.

\begin{figure}[H]
\begin{minipage}[t]{0.49\linewidth}
\centering
\includegraphics[width=0.8\linewidth]{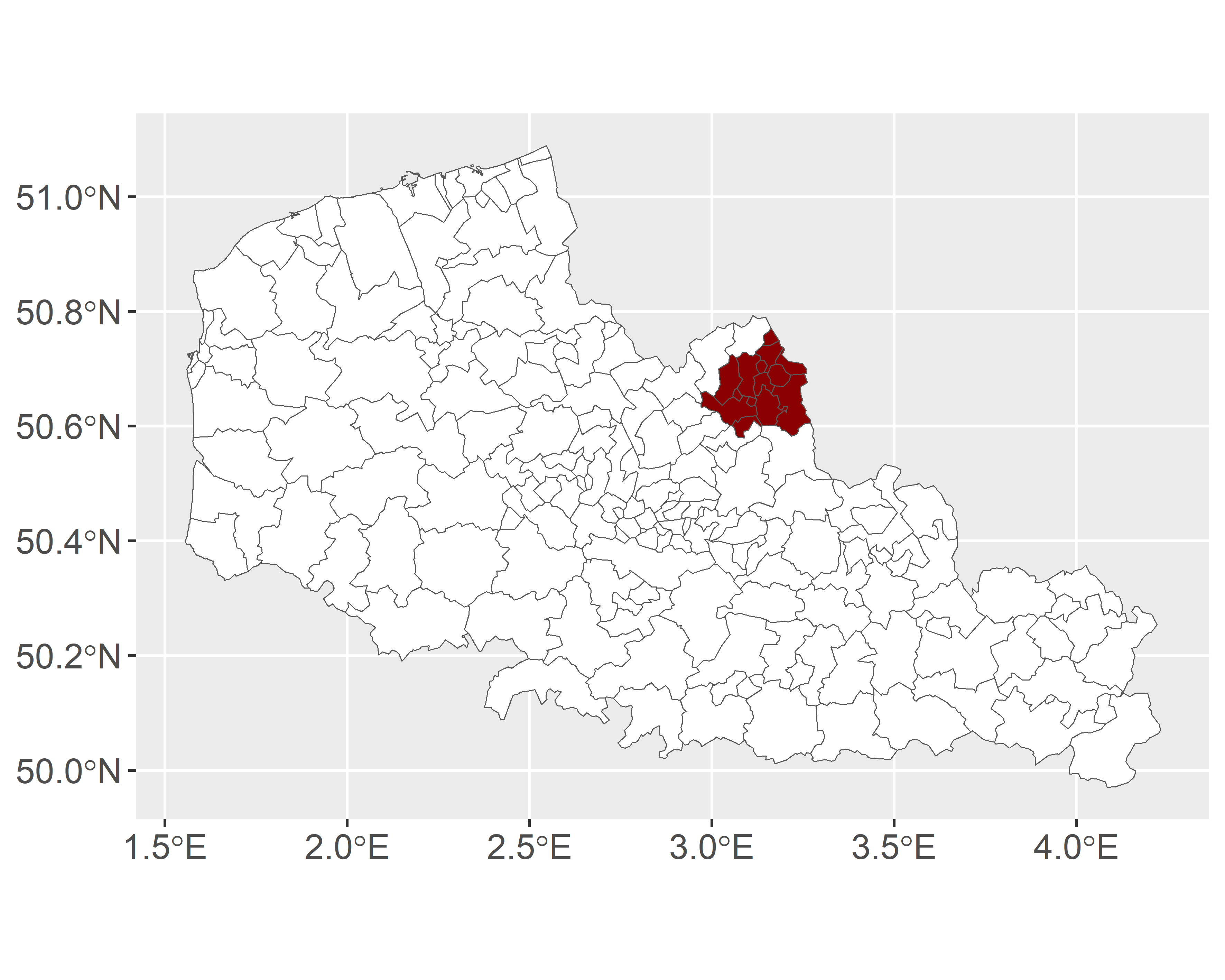}
\end{minipage}
\hfill
\begin{minipage}[t]{0.49\linewidth}
\centering
\includegraphics[width=0.8\linewidth]{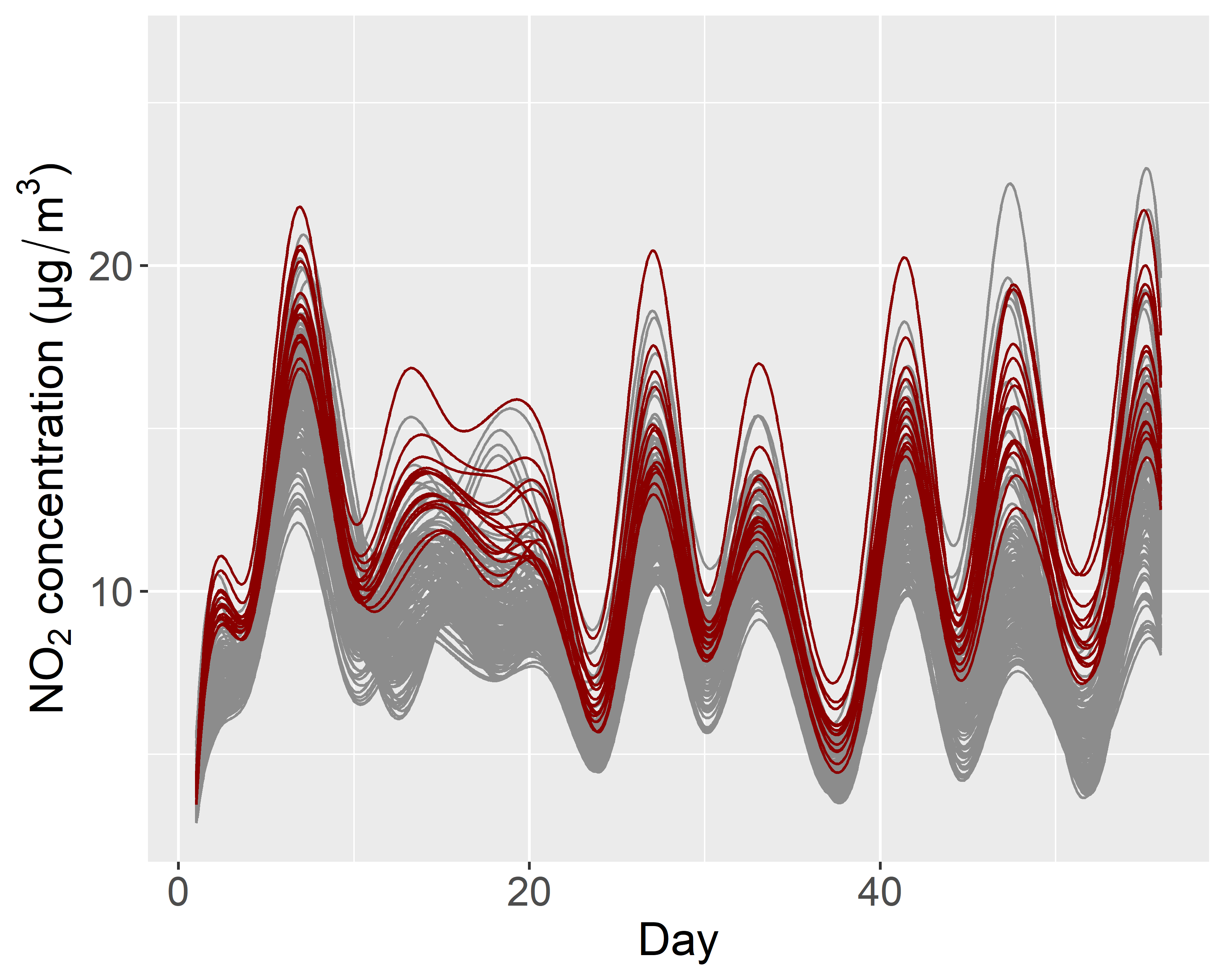}
\end{minipage}
\caption{Visualization of the detected cluster with the DFFSS method as well as the concentration curves of $\text{NO}_2$ (in $\mu g/m^3$) inside (in red) and outside (in grey) the cluster}
\label{fig:pollutioncluster}
\end{figure}
Figure \ref{fig:pollutioncluster} (left panel) displays the statistically significant cluster detected by the DFFSS (highlighted in red) in this air pollution dataset. The right panel compares the $\text{NO}_2$ concentration over the time within this cluster (in red) with that outside the cluster (in grey), revealing higher concentrations of $\text{NO}_2$ within the identified cluster. This information can be valuable for authorities in conducting local investigations and implementing policies to mitigate pollution. \\

FDA's significance has grown significantly owing to its relevance across diverse domains and the advancements in data collection technology (please refer to \cite{koner2023second} for further insights). The role of FDA in comprehending, analyzing, and harnessing these datasets is poised to expand further. FDA's utilization in burgeoning fields will influence the trajectory of data-driven innovation, decision-making, and issue resolution. Here is a glimpse of FDA's potential past and forthcoming contributions in these areas:

\begin{itemize}
\item \textit{Healthcare and Personalized Medicine}:
FDA can analyze patient data as continuous functions, allowing for personalized treatment plans based on individual health profiles. Real-time monitoring through wearables and FDA can aid in disease prediction and the optimization of treatment strategies.

\item \textit{Artificial Intelligence (AI) and Machine Learning (ML)}:
FDA provides a nuanced representation of data, improving AI and ML models' performance in various applications. In speech recognition, it enhances accuracy by capturing the continuous nature of speech signals.

\item \textit{Climate Science}:
In high-resolution climate data analysis, FDA identifies subtle patterns and trends, aiding in modeling, prediction, and mitigation strategies. Continuous data analysis contributes to precise climate projections and monitoring environmental changes.

\item \textit{Digital Marketing and User Behavior}:
In the digital realm, FDA uncovers intricate user behavior patterns, optimizing marketing, user experience, and product recommendations. It analyzes continuous data streams from digital platforms for deeper insights.

\item \textit{Brain-Computer Interfaces (BCIs)}:
FDA enhances BCIs by interpreting continuous brain activity data for prosthetics, neurorehabilitation, and cognitive augmentation. It enables precise control of assistive devices and cognitive enhancements.

\item \textit{Smart Cities}:
In smart cities, FDA optimizes urban planning, transportation systems, and energy consumption by analyzing continuous IoT and sensor data. It helps design sustainable and efficient cities through traffic analysis and energy usage trends.

\item \textit{Biotechnology and Synthetic Biology}:
In biotechnology, FDA models complex biological systems, facilitating the design of custom organisms and pharmaceuticals. It analyzes longitudinal data to engineer organisms for specific tasks.
\end{itemize}

\paragraph{Softwares:}
The Task View Functional Data Analysis on CRAN (\url{https://cran.r-project.org/web/views/FunctionalData.html}) lists the available \textsf{R} packages in the field of FDA, covering general functional data analysis, unsupervised learning (PCA, clustering, \dots), supervised learning (regression, classification), visualization and exploratory data analysis, registration, and alignment. \textsf{Python} and \textsf{MATLAB} also offer a few alternatives such as \textit{fdasrsf} and \textit{scikit-fda} (for \textsf{Python}), and \textit{fda} or \textit{fdasrvf} (for \textsf{MATLAB}).

\section{Conclusion}

Functional Data Analysis (FDA) understands and extracts meaningful insights from data that continuously vary over a continuum. While FDA may be particularly intriguing for those with a mathematical inclination, it invites everyone to explore the process of transforming numbers into valuable insights and offers a statistical approach that allows us to gain a deeper understanding of the world and actively contribute to shaping the future.

Indeed, over time, the methods and frequency of data collection will evolve, and computing and storage capacities will increase. The development of functional analysis methods is therefore essential, and their applications will improve decision-making in a variety of fields, providing biologists, economists, and policymakers with accurate information to make informed choices.

FDA is therefore a powerful field for understanding, analyzing, and using complex datasets. Thus, next time you see a graph, don't just see points and lines, but look for the continuous story it tells, the hidden patterns it holds, and the insights it offers. Remember what you just read: this is the realm of functional data analysis, where numbers transform into narratives waiting to be discovered.

\bibliographystyle{chicago}
\bibliography{bibliographie.bib}

\begin{thebibliography}{}

\bibitem[\protect\citeauthoryear{Berrendero, Justel, and Svarc}{Berrendero
  et~al.}{2011}]{berrendero2011principal}
Berrendero, J.~R., A.~Justel, and M.~Svarc (2011).
\newblock Principal components for multivariate functional data.
\newblock {\em Computational Statistics \& Data Analysis\/}~{\em 55\/}(9),
  2619--2634.

\bibitem[\protect\citeauthoryear{Bosq}{Bosq}{2000}]{bosq2000linear}
Bosq, D. (2000).
\newblock {\em Linear processes in function spaces: theory and applications},
  Volume 149.
\newblock Springer Science \& Business Media.

\bibitem[\protect\citeauthoryear{Cremona, Xu, Makova, Reimherr, Chiaromonte,
  and Madrigal}{Cremona et~al.}{2019}]{cremona2019functional}
Cremona, M.~A., H.~Xu, K.~D. Makova, M.~Reimherr, F.~Chiaromonte, and
  P.~Madrigal (2019).
\newblock Functional data analysis for computational biology.
\newblock {\em Bioinformatics\/}~{\em 35\/}(17), 3211.

\bibitem[\protect\citeauthoryear{Cuevas, Febrero, and Fraiman}{Cuevas
  et~al.}{2002}]{cuevas2002linear}
Cuevas, A., M.~Febrero, and R.~Fraiman (2002).
\newblock Linear functional regression: the case of fixed design and functional
  response.
\newblock {\em Canadian Journal of Statistics\/}~{\em 30\/}(2), 285--300.

\bibitem[\protect\citeauthoryear{Ferraty and Vieu}{Ferraty and
  Vieu}{2006}]{ferraty2006nonparametric}
Ferraty, F. and P.~Vieu (2006).
\newblock {\em Nonparametric functional data analysis}.
\newblock Springer.

\bibitem[\protect\citeauthoryear{Fr{\'e}vent, Ahmed, Dabo-Niang, and
  Genin}{Fr{\'e}vent et~al.}{2023}]{frevent2023investigating}
Fr{\'e}vent, C., M.-S. Ahmed, S.~Dabo-Niang, and M.~Genin (2023).
\newblock Investigating spatial scan statistics for multivariate functional
  data.
\newblock {\em Journal of the Royal Statistical Society Series C: Applied
  Statistics\/}~{\em 72\/}(2), 450--475.

\bibitem[\protect\citeauthoryear{Fr{\'e}vent, Ahmed, Marbac, and
  Genin}{Fr{\'e}vent et~al.}{2021}]{frevent2021detecting}
Fr{\'e}vent, C., M.-S. Ahmed, M.~Marbac, and M.~Genin (2021).
\newblock Detecting spatial clusters in functional data: New scan statistic
  approaches.
\newblock {\em Spatial Statistics\/}~{\em 46}, 100550.

\bibitem[\protect\citeauthoryear{Garc{\'\i}a, Garc{\'\i}a-R{\'o}denas, and
  G{\'o}mez}{Garc{\'\i}a et~al.}{2015}]{garcia2015k}
Garc{\'\i}a, M. L.~L., R.~Garc{\'\i}a-R{\'o}denas, and A.~G. G{\'o}mez (2015).
\newblock K-means algorithms for functional data.
\newblock {\em Neurocomputing\/}~{\em 151}, 231--245.

\bibitem[\protect\citeauthoryear{Hall, M{\"u}ller, and Wang}{Hall
  et~al.}{2006}]{10.1214/009053606000000272}
Hall, P., H.-G. M{\"u}ller, and J.-L. Wang (2006).
\newblock {Properties of principal component methods for functional and
  longitudinal data analysis}.
\newblock {\em The Annals of Statistics\/}~{\em 34\/}(3), 1493 -- 1517.

\bibitem[\protect\citeauthoryear{Hitchcock, Booth, and Casella}{Hitchcock
  et~al.}{2007}]{hitchcock2007effect}
Hitchcock, D.~B., J.~G. Booth, and G.~Casella (2007).
\newblock The effect of pre-smoothing functional data on cluster analysis.
\newblock {\em Journal of Statistical Computation and Simulation\/}~{\em
  77\/}(12), 1043--1055.

\bibitem[\protect\citeauthoryear{Horv{\'a}th and Kokoszka}{Horv{\'a}th and
  Kokoszka}{2012}]{horvath2012inference}
Horv{\'a}th, L. and P.~Kokoszka (2012).
\newblock {\em Inference for functional data with applications}, Volume 200.
\newblock Springer Science \& Business Media.

\bibitem[\protect\citeauthoryear{James}{James}{2002}]{james2002generalized}
James, G.~M. (2002).
\newblock Generalized linear models with functional predictors.
\newblock {\em Journal of the Royal Statistical Society Series B: Statistical
  Methodology\/}~{\em 64\/}(3), 411--432.

\bibitem[\protect\citeauthoryear{Karhunen}{Karhunen}{1947}]{karhunen1947under}
Karhunen, K. (1947).
\newblock Under lineare methoden in der wahr scheinlichkeitsrechnung.
\newblock {\em Annales Academiae Scientiarun Fennicae Series A1: Mathematia
  Physica\/}~{\em 47}.

\bibitem[\protect\citeauthoryear{Kleffe}{Kleffe}{1973}]{kleffe1973principal}
Kleffe, J. (1973).
\newblock Principal components of random variables with values in a seperable
  hilbert space.
\newblock {\em Mathematische Operationsforschung und Statistik\/}~{\em 4\/}(5),
  391--406.

\bibitem[\protect\citeauthoryear{Kokoszka and Reimherr}{Kokoszka and
  Reimherr}{2017}]{kokoszka2017introduction}
Kokoszka, P. and M.~Reimherr (2017).
\newblock {\em Introduction to functional data analysis}.
\newblock CRC press.

\bibitem[\protect\citeauthoryear{Koner and Staicu}{Koner and
  Staicu}{2023}]{koner2023second}
Koner, S. and A.-M. Staicu (2023).
\newblock Second-generation functional data.
\newblock {\em Annual Review of Statistics and Its Application\/}~{\em 10},
  547--572.

\bibitem[\protect\citeauthoryear{Kosambi}{Kosambi}{1943}]{kosambi2016decomp}
Kosambi, D.~D. (1943).
\newblock Statistics in function space.
\newblock {\em Journal of the Indian Mathematical Society\/}, 76--88.

\bibitem[\protect\citeauthoryear{Leng and M{\"u}ller}{Leng and
  M{\"u}ller}{2006}]{leng2006classification}
Leng, X. and H.-G. M{\"u}ller (2006).
\newblock Classification using functional data analysis for temporal gene
  expression data.
\newblock {\em Bioinformatics\/}~{\em 22\/}(1), 68--76.

\bibitem[\protect\citeauthoryear{Locantore, Marron, Simpson, Tripoli, Zhang,
  Cohen, Boente, Fraiman, Brumback, Croux, et~al.}{Locantore
  et~al.}{1999}]{locantore1999robust}
Locantore, N., J.~Marron, D.~Simpson, N.~Tripoli, J.~Zhang, K.~Cohen,
  G.~Boente, R.~Fraiman, B.~Brumback, C.~Croux, et~al. (1999).
\newblock Robust principal component analysis for functional data.
\newblock {\em Test\/}~{\em 8}, 1--73.

\bibitem[\protect\citeauthoryear{Loève}{Loève}{1945}]{loeve1945functions}
Loève, M. (1945).
\newblock Fonctions aléatoires du second ordre.
\newblock {\em Comptes Rendus Académie des Sciences\/}, 220--380.

\bibitem[\protect\citeauthoryear{Mercer}{Mercer}{1909}]{mercer1909functions}
Mercer, J. (1909).
\newblock Functions of positive and negative type and their connection with
  theory of integral equations.
\newblock {\em Philosophical transactions of the royal society of London.
  Series A\/}, 209(441--458):415–446.

\bibitem[\protect\citeauthoryear{Morris}{Morris}{2015}]{morris2015functional}
Morris, J.~S. (2015).
\newblock Functional regression.
\newblock {\em Annual Review of Statistics and Its Application\/}~{\em 2},
  321--359.

\bibitem[\protect\citeauthoryear{Palumbo, Centofanti, and Del~Re}{Palumbo
  et~al.}{2020}]{fabioetal20}
Palumbo, B., F.~Centofanti, and F.~Del~Re (2020).
\newblock Function-on-function regression for assessing production quality in
  industrial manufacturing.
\newblock {\em Quality and Reliability Engineering International\/}~{\em
  36\/}(8), 2738--2753.

\bibitem[\protect\citeauthoryear{Ramsay and Silverman}{Ramsay and
  Silverman}{2005}]{ramsay2005}
Ramsay, J.~O. and B.~W. Silverman (2005).
\newblock {\em Functional Data Analysis\/} (Second ed.).
\newblock Springer Series in Statistics. Springer.

\bibitem[\protect\citeauthoryear{Sangalli, Secchi, Vantini, and
  Vitelli}{Sangalli et~al.}{2010}]{sangalli2010k}
Sangalli, L.~M., P.~Secchi, S.~Vantini, and V.~Vitelli (2010).
\newblock K-mean alignment for curve clustering.
\newblock {\em Computational Statistics \& Data Analysis\/}~{\em 54\/}(5),
  1219--1233.

\bibitem[\protect\citeauthoryear{Shang}{Shang}{2014}]{shang2014survey}
Shang, H.~L. (2014).
\newblock A survey of functional principal component analysis.
\newblock {\em AStA Advances in Statistical Analysis\/}~{\em 98}, 121--142.

\bibitem[\protect\citeauthoryear{Silverman and Ramsay}{Silverman and
  Ramsay}{2002}]{silverman2002applied}
Silverman, B. and J.~Ramsay (2002).
\newblock Applied functional data analysis: methods and case studies.

\bibitem[\protect\citeauthoryear{Wang, Chiou, and M{\"u}ller}{Wang
  et~al.}{2016}]{wang2016functional}
Wang, J.-L., J.-M. Chiou, and H.-G. M{\"u}ller (2016).
\newblock Functional data analysis.
\newblock {\em Annual Review of Statistics and its application\/}~{\em 3},
  257--295.

\bibitem[\protect\citeauthoryear{Xiao}{Xiao}{2020}]{10.3150/20-BEJ1209}
Xiao, L. (2020).
\newblock {Asymptotic properties of penalized splines for functional data}.
\newblock {\em Bernoulli\/}~{\em 26\/}(4), 2847 -- 2875.

\bibitem[\protect\citeauthoryear{Yao, M{\"u}ller, and Wang}{Yao
  et~al.}{2005}]{yao2005functional}
Yao, F., H.-G. M{\"u}ller, and J.-L. Wang (2005).
\newblock Functional data analysis for sparse longitudinal data.
\newblock {\em Journal of the American statistical association\/}~{\em
  100\/}(470), 577--590.

\bibitem[\protect\citeauthoryear{Zhang and Parnell}{Zhang and
  Parnell}{2023}]{zhang2023review}
Zhang, M. and A.~Parnell (2023).
\newblock Review of clustering methods for functional data.
\newblock {\em ACM Transactions on Knowledge Discovery from Data\/}~{\em
  17\/}(7), 1--34.

\bibitem[\protect\citeauthoryear{Zhang and Wang}{Zhang and
  Wang}{2016}]{zhangwang}
Zhang, X. and J.-L. Wang (2016).
\newblock {From sparse to dense functional data and beyond}.
\newblock {\em The Annals of Statistics\/}~{\em 44\/}(5), 2281 -- 2321.

\end{thebibliography}

\end{document}